\documentclass{article}

\usepackage{times}
\usepackage{graphics}
\usepackage{amssymb}
\usepackage{natbib}

\newcommand{\dif}{\mathrm{d}}

\newcommand{\Beta}{\mathrm{Beta}}

\newcommand{\Gma}{\mathrm{Gamma}}

\newcommand{\pr}{\Pr}

\begin{document}

\title{A closed-form approximation for the median of the beta distribution}
\author{Jouni Kerman}
\date{November 1, 2011}
\maketitle


\begin{abstract}
A simple closed-form approximation for the
median of the beta distribution $\Beta(a, b)$
is introduced:
$(a-1/3)/(a+b-2/3)$ for $(a,b)$ both larger than $1$
has a relative error of less than 4\%,
rapidly decreasing to zero as both shape parameters increase.
\end{abstract}


\section{Introduction}
\label{s1}

Consider 
the the beta distribution
$\Beta(a,b)$, with the density function,
$$
  \frac{\Gamma(a+b)}{\Gamma(a)\Gamma(b)}
  \theta^{a-1}(1-\theta)^{b-1}.
$$
The mean of $\Beta(a, b)$ is readily obtained
by the formula $a/(a+b)$, but
there is no general closed formula for the median.
The median function, here denoted by $m(a,b)$,
is the function that satisfies,
$$
  \frac{\Gamma(a+b)}{\Gamma(a)\Gamma(b)}
  \int_0^{m(a,b)}\theta^{a-1}(1-\theta)^{b-1} \dif\theta
  =
  \frac{1}{2}.
$$

The relationship $m(a,b)=1-m(b, a)$ holds.
Only for the special cases $a=1$ or $b=1$ we may obtain 
an exact formula: $m(a, 1)=2^{-1/a}$
and $m(1, b)=1-2^{-1/b}$.
Moreover, when $a=b$, the median is exactly $1/2$.

There has been much literature about the 
incomplete beta function and 
its inverse
(see e.g. \citet{Dutka:1981} for a review).
The focus in literature has been
on finding accurate numerical results,
but a simple and practical 
approximation that is easy to compute 
has not been found.

\begin{figure}[t]
\begin{center}
\scalebox{0.80}{\includegraphics{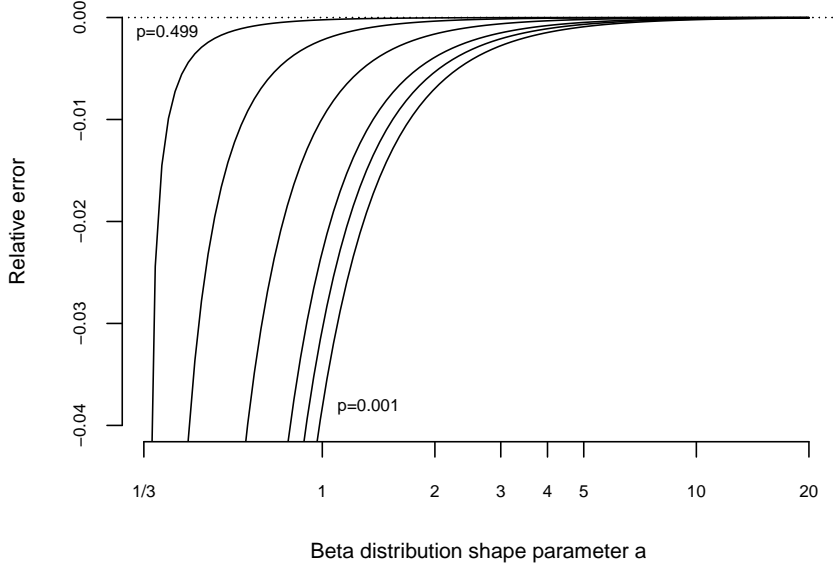}}
\caption{\label{fig-betaerrors}
Relative errors of the approximation
$(a-1/3)/(a+b-2/3)$ of the 
median of the $\Beta(a, b)$ distribution,
compared with the numerically computed value
for several fixed $p=a/(a+b)<1/2$.
The horizontal axis shows the shape parameter $a$
on logarithmic scale.
From left to right, 
$p=0.499$, 0.49, 0.45, 0.35, 0.25, and 0.001.
}
\end{center}
\end{figure}
\section{A new closed-form approximation for the median}

Trivial bounds for the median can be derived
\citep{Payton:1989}, which are 
a consequence of the more general
mode-median-mean inequality 
\citep{Groeneveld:Meeden:1977}.
In the case of the beta distribution with
$1<a<b$,
the median is bounded by the 
mode $(a-1)/(a+b-2)$ and the mean $a/(a+b)$:
$$
  \frac{a-1}{a+b-2} 
  \le
  m(a,b)
  \le
  \frac{a}{a+b}.
$$
For $a\le1$ the formula for the mode does not hold
as there is no mode.
If $1<b<a$, the order of the inequality is reversed.
Equality holds if and only if $a=b$;
in this case the mean, median, and mode are all equal to $1/2$.

This inequality shows that if the mean is kept fixed
at some $p$, 
and one of the shape parameters is increased, say $a$, 
then the median is sandwiched between
$p(a-1)/(a-2p)$ and $p$,
hence the median tends to $p$.


From the formulas for the mode and mean,
it can be conjectured that 
the median $m(a,b)$ could be approximated by
$m(a,b;d)=(a-d)/(a+b-2d)$ for some $d\in(0,1)$,
as this form would satisfy the above inequality
while agreeing with the symmetry requirement, 
that is, $m(a,b;d)=1-m(b,a;d)$.

\begin{figure}[t]
\begin{center}
\scalebox{0.80}{\includegraphics{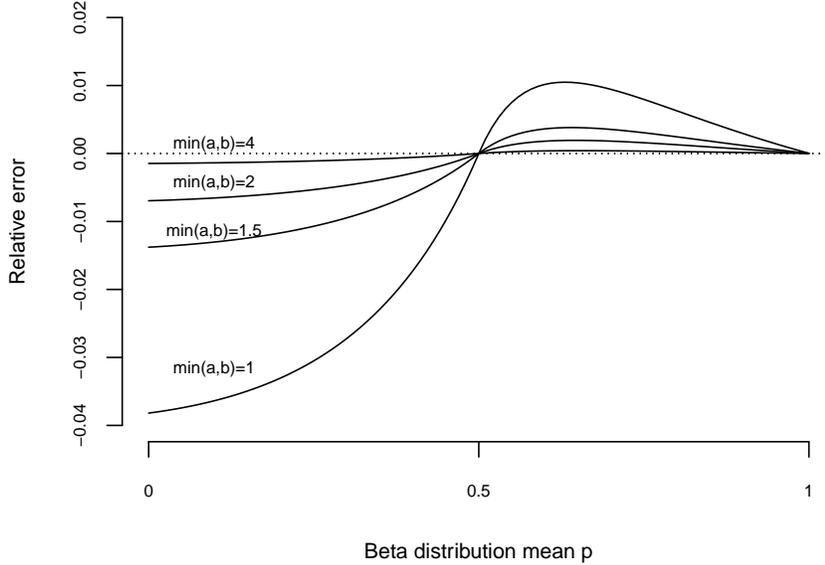}}
\caption{\label{fig-betaerrp}
Relative errors of the approximation
$(a-1/3)/(a+b-2/3)$ of the 
median of the $\Beta(a, b)$ distribution
over the whole range of possible distribution means
$p=a/(a+b)$.
The smaller of the shape parameters is fixed,
i.e. for $p\le 0.5$, 
the median is computed for $\Beta(a, a(1-p)/p)$
and for $p>0.5$, 
the median is computed for $\Beta(bp/(1-p), b)$.
}
\end{center}
\end{figure}

Since a $\Beta(a,b)$ variate can be expressed as 
the ratio $\gamma_1/(\gamma_1+\gamma_2)$ where
$\gamma_1\sim\Gma(a)$ and 
$\gamma_2\sim\Gma(b)$ (both with unit scale),
it is useful to have a look at the median 
of the gamma distribution. 
\citet{Berg:Pedersen:2006} studied the median 
function of the unit-scale 
gamma distribution median function, denoted here by $M(a)$,
for any shape parameter $a>0$,
and obtained
$M(a) = a - 1/3 + o(1)$,
rapidly approaching $a-1/3$ as $a$ increases.
It can therefore be conjectured that
the distribution median may be approximated by,
\begin{equation}\label{eq-beta}
 m(a, b) \approx
 m(a, b; 1/3)
 = 
 \frac{a-1/3}{(a-1/3)+(b-1/3)}
 =
 \frac{a-1/3}{a+b-2/3}.
\end{equation}

Figure (\ref{fig-betaerrors}) 
shows that 
this approximation indeed appears to approach 
the numerically computed median asymptotically 
for all distribution means $p=a/(a+b)$ as the
(smaller) shape parameter $a\to\infty$. 
For $a\ge1$, the relative error is less than 4\%,
and for $a\ge2$ this is already less than 1\%.

\begin{figure}[t]
\begin{center}
\scalebox{0.80}{\includegraphics{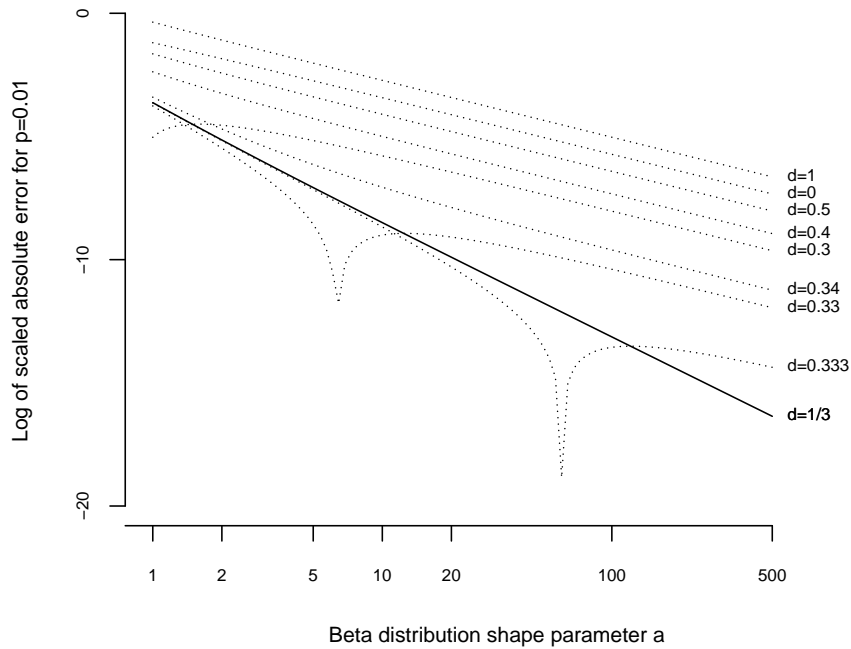}}
\caption{\label{fig-betadisterr}
Logarithm of the scaled absolute error (distance)
$\log(|m(a,b;d)-m(a,b)|/p)$,
computed for a fixed distribution mean $p=0.01$ and various 
$d$. The approximate median of the
$\Beta(a,b)$ distribution
is defined as $m(a,b;d)=(a-d)/(a+b-2d)$.
Due to scaling of the error, 
the graph and its scale will not essentially change even if 
the error is computed for other values of $p<0.5$.
The approximation $m(a,b;1/3)$ performs the 
most consistently, 
attaining the lowest absolute error eventually as the
precision of the distribution increases.
}
\end{center}
\end{figure}
Figure (\ref{fig-betaerrp}) shows the relative
error over all possible distribution means $p=a/(a+b)$,
as the smallest of the two shape parameters varies from
$1$ to $4$. This illustrates how the relative error
tends uniformly to zero over all $p$ as the shape parameters 
increase. 
The figure also shows that
the formula consistently either underestimates
or overestimates the median depending on whether
$p<0.5$ or $p>0.5$.

However, the function
$m(a,b;d)$
approximates the median fairly accurately
if some other $d$ close to $1/3$ (say $d=0.3$) is chosen.
Figure (\ref{fig-betadisterr}) displays 
curves of the logarithm of the absolute
difference from the numerically computed
median for a fixed $p=0.01$, as the shape parameter
$a$ increases. 
The absolute difference
has been scaled by $p$ before taking the logarithm:
due to this scaling, 
the error stays approximately constant as $p$ decreases
so the picture and its scale will not essentially change even if 
the error is computed for other values of $p<0.5$.
The figure shows that although some approximations such as
$d=0.3$ has a lower absolute error for some $a$, 
the error of $m(a, b; 1/3)$ tends to be lower in the long run,
and moreover performs more consistently 
by decreasing at the same rate on the logarithmic scale.
In practical applications, $d=0.333$ should be a sufficiently
good approximation of $d=1/3$.

\begin{figure}[t]
\begin{center}
\scalebox{0.7}{\includegraphics{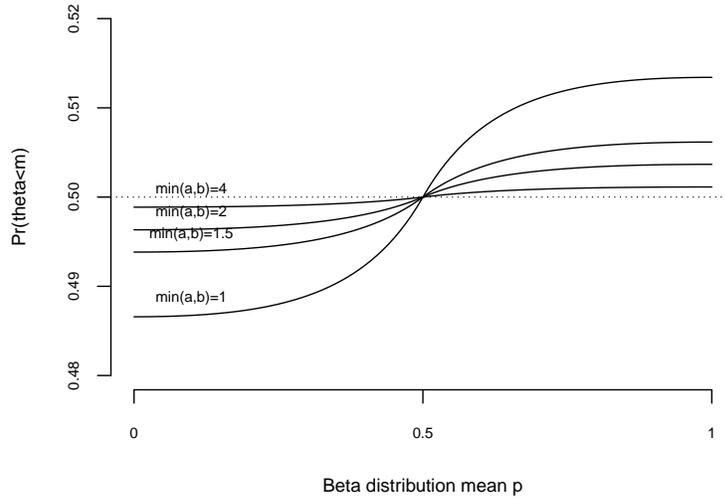}}
\caption{\label{fig-betatail}
Tail probabilities $\pr(\theta<m)$
of the $\Beta(a, b)$ distribution
when $m=(a-1/3)/(a+b-2/3)$.
As the smaller of the two shape
parameters increases, the tail probability
tends rapidly and uniformly to $0.5$.
}
\end{center}
\end{figure}
Another measure of the accuracy is the
tail probability
$\pr(\theta \le m(a,b;1/3))$ of a $\Beta(a, b)$ variate $\theta$:
good approximators of the 
median should yield probabilities close to $1/2$.
Figure (\ref{fig-betatail}) shows that
as long as the smallest of the shape parameters 
is at least 1,
the tail probability is bound between $0.4865$ and $0.5135$.
As the shape parameters increase, the 
probability tends 
rapidly and uniformly to $0.5$.

Finally, let us have a look at a
well-known paper that provides further
support for the uniqueness of $m(a,b;1/3)$.
\citet{Peizer:Pratt:1968} and \citet{Pratt:1968}
provide approximations for 
the probability function $\pr(\theta\le x)$
of a $\Beta(a,b)$ variate $\theta$.
Although they do not provide a formula
for the inverse, it is 
the probability function at the approximate median. 
According to \citet{Peizer:Pratt:1968},
$\pr(\theta\le x)$ 
is well approximated by
$\Phi(z(a,b; x))$ where $\Phi$ is the
standard normal probability function,
and $z$ is a function of the shape parameters and 
the quantile $x$. 
Consider $m=m(a,b;d)$:
$z(a,b;m)$ should be close to zero and at least
tend to zero fast as $a$ and $b$ increase. 
Now assume that $p$ is fixed, $a$ varies and $b=a(1-p)/p$. 
The function $z(a, b; m)$ equals, rewritten with
the notation in this paper,
\begin{equation}\label{eq-peizer-beta}
   \sqrt{p}\frac{1-2m}{(a-p)^{1/2}}\left(
    1/3-d
    -
   \frac{0.02p}{a}\left[
   \frac{1}{2} + \frac{1-dp/a}{p(1-p)}
   \right]
   \right)
   \left(\frac{1+f(a,p;d)}{m(1-m)}\right)^{1/2},
\end{equation}
where the function $f(a,p;d)$ tends to zero
as $a$ increases,
being exactly zero only when $d=1/2$ or $m=1/2$.
It is evident that for the fastest convergence
rate to zero, one should choose $d=1/3$.
This is of the order $O(a^{-3/2})$;
if $d\ne 1/3$, 
for example if we choose the mean $p$
as the approximation of the median ($d=0$),
the rate is at most $O(a^{-1/2})$.


\begin{thebibliography}{6}
\expandafter\ifx\csname natexlab\endcsname\relax\def\natexlab#1{#1}\fi
\expandafter\ifx\csname url\endcsname\relax
  \def\url#1{\texttt{#1}}\fi
\expandafter\ifx\csname urlprefix\endcsname\relax\def\urlprefix{URL }\fi
\providecommand{\eprint}[2][]{\url{#2}}
\providecommand{\bibinfo}[2]{#2}
\ifx\xfnm\relax \def\xfnm[#1]{\unskip,\space#1}\fi
\bibitem[{Berg and Pedersen(2006)}]{Berg:Pedersen:2006}
\bibinfo{author}{Berg, C.}, \bibinfo{author}{Pedersen, H.L.},
  \bibinfo{year}{2006}.
\newblock \bibinfo{title}{The {C}hen-{R}ubin conjecture in a continuous
  setting}.
\newblock \bibinfo{journal}{Methods and Applications of Analysis}
  \bibinfo{volume}{13}, \bibinfo{pages}{63--88}.
\bibitem[{Dutka(1981)}]{Dutka:1981}
\bibinfo{author}{Dutka, J.}, \bibinfo{year}{1981}.
\newblock \bibinfo{title}{The incomplete beta function -- a historical
  profile}.
\newblock \bibinfo{journal}{Archive for history of exact sciences}
  \bibinfo{volume}{24}, \bibinfo{pages}{11--29}.
\bibitem[{Groeneveld and Meeden(1977)}]{Groeneveld:Meeden:1977}
\bibinfo{author}{Groeneveld, R.A.}, \bibinfo{author}{Meeden, G.},
  \bibinfo{year}{1977}.
\newblock \bibinfo{title}{The mode, median, and mean inequality}.
\newblock \bibinfo{journal}{The American Statistician} \bibinfo{volume}{31},
  \bibinfo{pages}{120--121}.
\bibitem[{Payton et~al.(1989)Payton, Young and Young}]{Payton:1989}
\bibinfo{author}{Payton, M.}, \bibinfo{author}{Young, L.},
  \bibinfo{author}{Young, J.}, \bibinfo{year}{1989}.
\newblock \bibinfo{title}{Bounds for the difference between median and mean of
  beta and negative binomial distributions}.
\newblock \bibinfo{journal}{Metrika} \bibinfo{volume}{36},
  \bibinfo{pages}{347--354}.
\bibitem[{Peizer and Pratt(1968)}]{Peizer:Pratt:1968}
\bibinfo{author}{Peizer, D.B.}, \bibinfo{author}{Pratt, J.W.},
  \bibinfo{year}{1968}.
\newblock \bibinfo{title}{A normal approximation for binomial, {F}, beta, and
  other common, related tail probabilities, {I}}.
\newblock \bibinfo{journal}{Journal of the American Statistical Association}
  \bibinfo{volume}{63}, \bibinfo{pages}{1416--1456}.
\bibitem[{Pratt(1968)}]{Pratt:1968}
\bibinfo{author}{Pratt, J.W.}, \bibinfo{year}{1968}.
\newblock \bibinfo{title}{A normal approximation for binomial, {F}, beta, and
  other common, related tail probabilities, {II}}.
\newblock \bibinfo{journal}{Journal of the American Statistical Association}
  \bibinfo{volume}{63}, \bibinfo{pages}{1457--1483}.

\end{thebibliography}
\end{document}